\documentclass[11pt]{article}
\usepackage{geometry} 
\usepackage{mathptmx}       
\usepackage{helvet}         
\usepackage{courier}        
\usepackage{type1cm}        
\usepackage{graphicx}        
\usepackage{multicol}        
\usepackage[bottom]{footmisc}
\usepackage{subfigure}
\usepackage{amsmath}
\usepackage{amssymb}
\usepackage{hyperref}
\usepackage{cleveref}

\title{Biot's parameters estimation in ultrasound propagation through cancellous bone}
\author{Miguel Angel Moreles \and Joaquin Pe\~{n}a \and Jose Angel Neria}
\date{} 

\begin{document}

\maketitle
\begin{center}
Centro de Investigaci\'{o}n en Matem\'{a}ticas\\
Jalisco s/n, Valenciana\\
Guanajuato, GTO 36240,Mexico\\
email:moreles,jose.neria,joaquin@cimat.mx
\end{center}
\begin{abstract}
Of interest is the characterization of a cancellous bone immersed in an acoustic fluid. The bone is placed between an ultrasonic point source and a receiver.  Cancellous bone is regarded as a porous medium saturated with fluid according to Biot's theory. This model is coupled with the fluid in an open pore configuration and solved by means of the Finite Volume Method. Characterization is posed as a Bayesian parameter estimation problem in Biot's model given pressure data collected at the receiver.  As a first step we present numerical results in 2D for signal recovery.  It is shown that as point estimators, the Conditional Mean outperforms the classical PDE-constrained minimization solution.
\end{abstract}
\tableofcontents

\section{Introduction}

Analysis of initial boundary value problems (IBVPs) usually consider well-posed problems, that is, problems where uniqueness and existence of the solution, as well as continuous dependence on the input data can be establish. Problems which lack any of these properties are ill-posed or inverse problems \cite{tarantola}. For example, problems arising in geophysics and medicine concern the determination of properties of some inaccessible region. The problem of interest in this work is of this sort. The properties to estimate are parameters of a saturated porous medium immersed in a fluid, a so called Biot's medium. The parameters are to be recovered from an ultrasound noisy signal. 

The Biot's medium of concern is a medulla saturated cancellous bone. The bone is placed between an acoustic source and receiver. A potential application of the estimated parameters is as an aid in the diagnostic of osteoporosis.

Research on the problem is very active. In Buchanan et al [10, 11, 12] the problem of inversion of parameters for a two-dimensional sample of trabecular bone is considered in a low frequency range ($f <100$ KHz). In these works the recovered parameters are $\phi$ (porosity), $\alpha$ (solid tortuosity), $K_b$ (bulk modulus of the porous skeletal frame) and $N$ (solid shear modulus). In [9] trabecular bone samples are characterized by solving the inverse problem using experimentally acquired signals. In this case the direct problem is solved by a modified Biot's model, for the case a one-dimensional block of trabecular bone saturated with water. The recovered parameters are $\phi$, $\alpha$, $\nu_b$ (Poisson ratio of the porous skeletal frame) and $E_b$ (Young's modulus of the porous skeletal frame).

In these cases, the inversion problem is posed as minimization problem in the least squares sense. So the maximum likelihood estimator is obtained. In this setting, previous information such as physically acceptable ranges, results of other experiments, etc., may not be incorporated in a natural way. On the other hand, as pointed out in Sebaa et al. [9], the solution of the inverse problem for all model parameters using
only data from the transmitted signal is difficult, if not impossible.  This in part due to the high computational cost of the
optimization of the objective function and partly because more experimental data is needed to obtain a unique solution.

In the present work we follow an alternative approach. We pose the problem as one of Bayesian estimation. We show that in this approach, it is possible to estimate the parameters involved in the Biot's model, and previous information can be incorporated. 

\section{The Biot's model and the parameter estimation problem}

\subsection{A clinical motivation}

Noninvasive techniques for assessing bone fracture risk, as well as bone fragility are of current interest. Bone can be characterized in two types, cancellous (spongy or trabecular) and cortical. There is an ongoing discussion on trabecular changes due to osteoporosis, in particular, thinning of the trabeculae. Consequently, early detection of these changes is a potential aid for diagnostics of osteoporosis. A promising noninvasive technique is by ultrasound propagation trough cancellous bone, see Zebaa et al \cite{sebaaparam}, and references therein. 

\begin{figure}[ht]
\centering
\includegraphics[width=290pt]{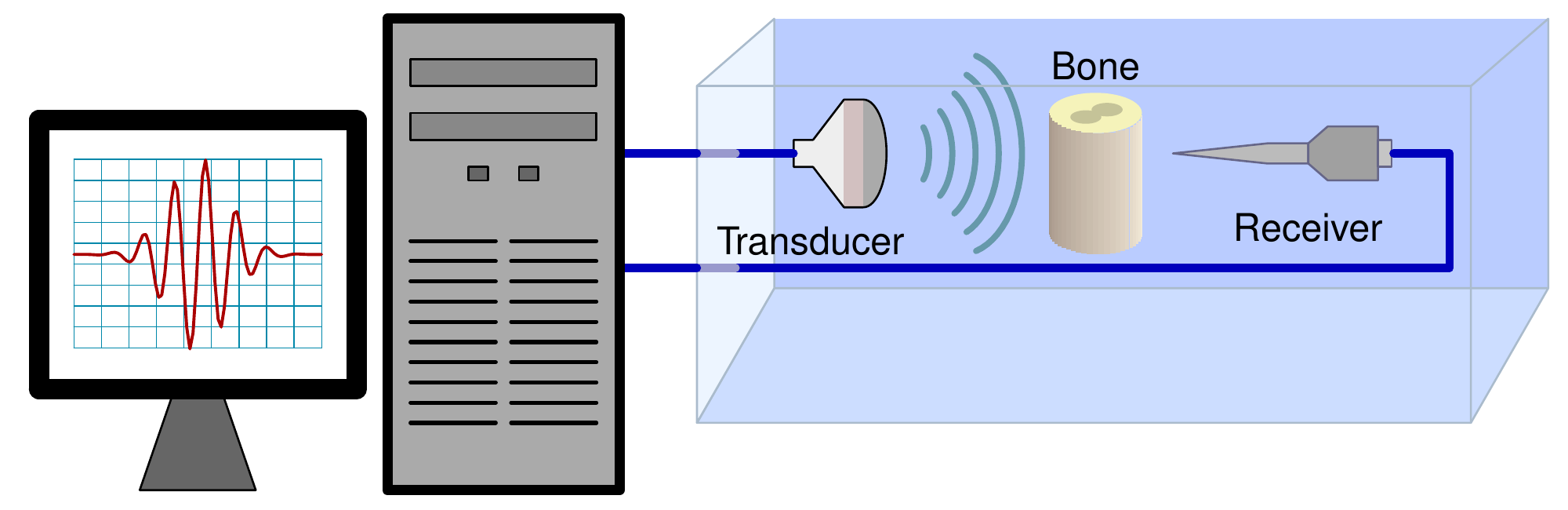}
\caption{Schematic showing the environment to study samples of cancellous bone using 
    ultrasonic axial transmission.}
\label{fig:through_trans}
\end{figure}

In accord with the \emph{axial transmission} (AT) technique, Lowet \& Van der \cite{LowetVander}, the configuration is a cancellous bone placed between an acoustic source and receiver.  A schematic is presented in Figure \ref{fig:through_trans}. 

An ultrasonic pulse is emitted from the transducer then propagated through the bone. The problem of interest, is to determine physical characteristics of the bone given a noisy signal collected at the receiver. 

\subsection{Biot's model for a fluid saturated porous solid}

The cancellous bone is modelled as a Biot medium, that is, a fluid saturated porous solid.  Biot's theory yields the following governing equations.

\begin{center}
\begin{subequations}\label{eqsBiot_perd_1}
\begin{equation}
\rho_{11} \frac{\partial^2  \vec{ U} ^s}{\partial t^2}        +\rho_{12}\frac{\partial^2  \vec{ U} ^f}{\partial t^2} =  \nabla \cdot  \sigma  - b\frac{\partial}{\partial t} \left(  \vec{ U}^s - \vec{ U}^f  \right)  ,
\end{equation} 
\begin{equation}
 \rho_{12} \frac{\partial^2  \vec{ U} ^s}{\partial t^2}  +   \rho_{22}\frac{\partial^2  \vec{ U} ^f}{\partial t^2}  = \nabla s  +   b\frac{\partial}{\partial t} \left(  \vec{ U}^s - \vec{ U}^f  \right)  ,
\end{equation} 
\begin{equation}
\sigma = [(P-2N)e+Q\varepsilon]I+2N\bar{e},
\end{equation}
\begin{equation}
s= Qe+R\varepsilon,
\end{equation}
\end{subequations}
\end{center}
where $\sigma$ and $s$ represent the forces acting  on  the  solid and  fluid
portions of each  side of an unit cube of the Biot medium,  respectively, 
$\vec{ U}^s$ and $\vec{ U}^f$ are solid and fluid displacements, and  

\begin{eqnarray*}
e             &=& \nabla\cdot \vec{ U}^s,\\
\varepsilon   &=& \nabla\cdot \vec{ U}^f,\\
\bar{e}_{i,j} &=& \frac{2-\delta_{i,j}}{2} \left( 
        \frac{\partial \vec{U}^s_i }{\partial x_j} + \frac{\partial \vec{U}^s_j}{\partial x_i} \right),\\
 I_{ij}       &=& \delta_{ij} = \begin{cases}
        0, \; & \text{if }\; i\neq j, \\
        1, \; & \text{if}\;  i=j.
                        \end{cases} \\
\end{eqnarray*}
\noindent
Also $P$, $Q$, $R$ are generalized elastic constants given by
\begin{subequations}\label{eq:generalized_ctt}
\begin{equation}
P=\frac{(1-\phi)\left( 1-\phi-\frac{K_b}{K_s}  \right) K_s + \phi \frac{K_f}{K_f}K_b  }{\Delta} + \frac{4N}{3},
\end{equation}
\begin{equation}
Q=\frac{ \left( 1-\phi-\frac{K_b}{K_s}  \right)\phi K_s  }{\Delta},
\end{equation}
\begin{equation}
R=\frac{\phi^2K_s}{\Delta},
\end{equation}
\begin{equation}
\Delta = 1-\phi- \frac{K_b}{K_s}+\phi \frac{K_s}{f}.
\end{equation}
\end{subequations}
The measurable quantities in these expressions are $\phi$ (porosity), $K_f$ (bulk modulus of the pore fluid), $K_s$ (bulk modulus of elastic solid) and $K_b$ (bulk modulus of the porous skeletal frame).  $N$ is the solid shear modulus.

The remaining parameters are the mass coupling coefficients, namely
\begin{subequations}\label{eq:mass_coef}
\begin{equation}
\rho_{11}+\rho_{12}=(1-\phi) \rho_\mathrm{s},
\end{equation}
\begin{equation}
\rho_{22}+\rho_{12}=\phi \rho_\mathrm{f},
\end{equation}
\begin{equation}
\rho_{12} = -(\alpha -1)\phi\rho_\mathrm{f}.
\end{equation}
\end{subequations}
where $\rho_\mathrm{s}$, $\rho_ \mathrm{f}$ are the solid and fluid densities, $\alpha$ is the solid tortuosity and $b$ is a parameter depending on the frequency of the incident wave and accounts for energy losses in the solid-fluid structure.

\subsection{Acoustic fluid}

\begin{figure}
\centering
\includegraphics[width=150pt]{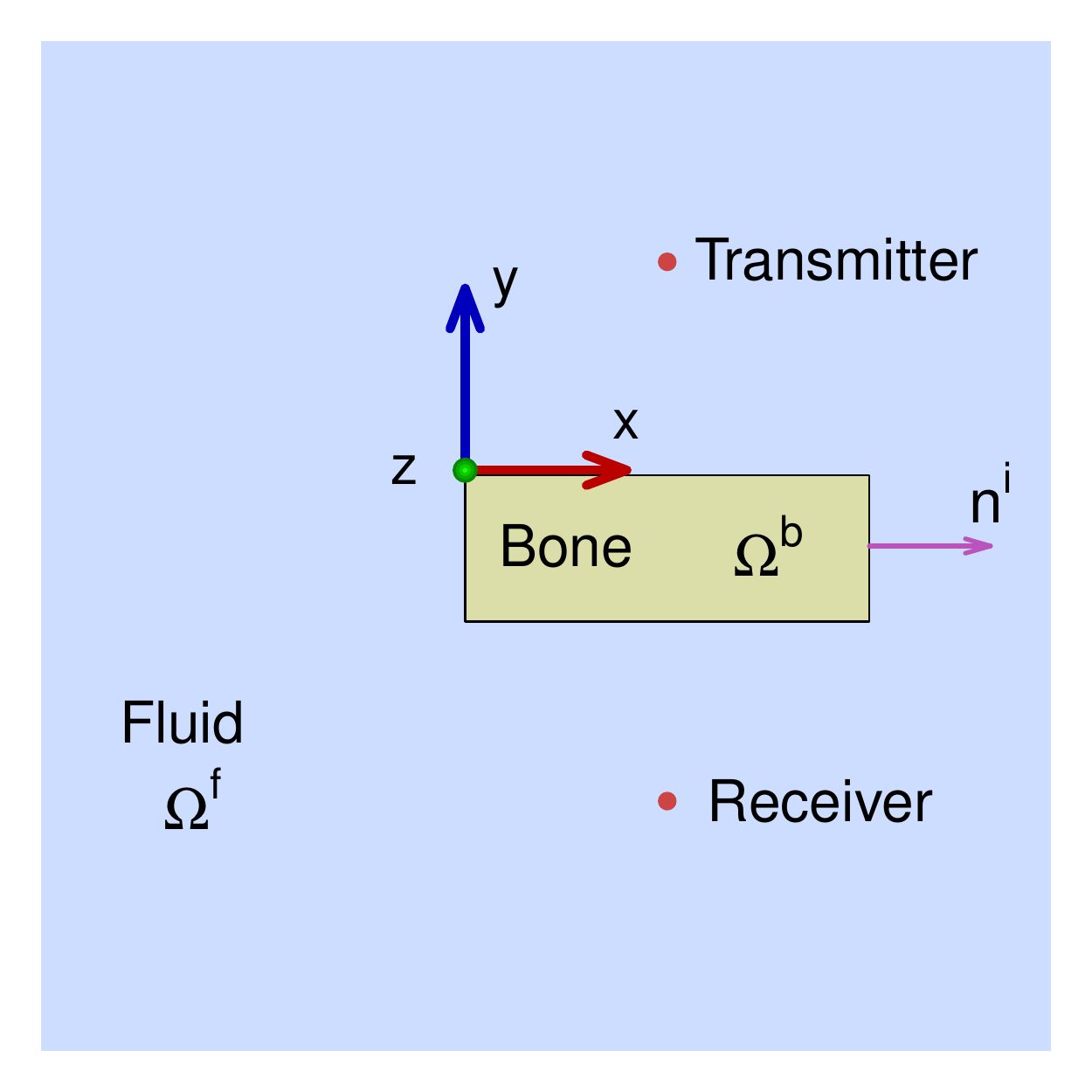}
\caption{Basic configuration of the elements in the domain.}
\label{fig:conf_directo}
\end{figure}

It is assumed that the Biot medium $\Omega^b$ (cancellous bone) is immersed in a fluid
as shown in Figure  \ref{fig:conf_directo}. 

The fluid within  $\Omega^f$ is acoustic. Its density and speed are $\rho$, $c$, respectively. Consequently, in terms of the pressure $P(\textbf{x}  \text{,} t)$ we have 

\begin{equation}\label{eq:scalar_wave}
\frac{1}{c^2}\frac{\partial^2  P}{\partial t^2} - \nabla^2 P = \frac{\partial Q}{\partial t}, 
\quad  \forall \vec{x}\ \in \ \Omega^f,
\end{equation}
where $Q(\textbf{x} \text{,}  t)$ is the point source density located at $\textbf{x}^s$ and given by
\begin{equation}
\frac{\partial Q}{\partial t} = \rho F(t)\delta(x-x^s)\delta(y-y^s),
\end{equation}
where $F(t)$ is a scalar real function and $\delta(\cdot)$ is the Dirac's delta function.

In $\Omega^f$, the velocity vector $\vec{v}(\textbf{x},t)$, is related with the pressure gradient by means of the Euler's equation 
\begin{equation}\label{eq:euler_pv}
\rho \frac{\partial \vec{v}}{\partial t}+\nabla P=\vec{0}, \quad  \forall \vec{x}\ \in \ \Omega^f.
\end{equation}

\subsection{Initial and boundary conditions}

According to Figure \ref{fig:conf_directo},  $\Omega^f$  is the domain occupied by the fluid whereas the fluid saturated porous medium is $\Omega^b$.  Null Dirichlet boundary conditions are prescribed in the outer boundary,

\begin{equation}\label{P_Bd}
P=0, \quad  \forall \vec{x}\ \in \ \partial\Omega^f.
\end{equation}

The configuration is chosen so that at the receiver, the waves propagating in the fluid are not affected.

As derived in Lovera \cite{LoveraBd}, in the Biot medium-fluid interface $\partial\Omega^b$,  the boundary conditions are  

\begin{equation}\label{eq:boundary_conds_new}
\left. \begin{array}{rcl}
  s & = & -\phi P \\
  \sigma\, \vec{n}^i & = & -(1- \phi)  P\,\vec{n}^i 
\end{array} \right\}
\quad \forall\; \vec{x}	 \  \in \ \partial\Omega^b,
\end{equation}
where $\vec{n}^i$ is the normal unit vector to the interface $\partial\Omega^b$ pointing from $\Omega^b$ towards the fluid. 

The system starts at rest. Consequently, zero initial conditions are added to the PDE system to have a well posed Initial Boundary Value Problem (IBVP). Namely

\begin{equation}\label{P_Ini}
P(\vec{x},0)=0, \quad  \forall \vec{x}\ \in \ \Omega^f.
\end{equation}

\begin{equation}\label{U_Ini}
\vec{ U}^s(\vec{x},0)=\vec{0}, \quad  \vec{ U}^f(\vec{x},0)=\vec{0}, \quad  \forall \vec{x}\ \in \ \Omega^s.
\end{equation}

\subsection{Numerical solution}
In parameter estimation problems, the solution of the so called forward map, in this case involving the solution of the IBVP, is taken for granted. We have made an implementation of the finite volume method \cite{mazumder}. Our choice is based on the easy handling of the 
 boundary conditions at the interface. Let us illustrate this fact. 
 
\begin{figure}
\centering
\includegraphics[width=150pt]{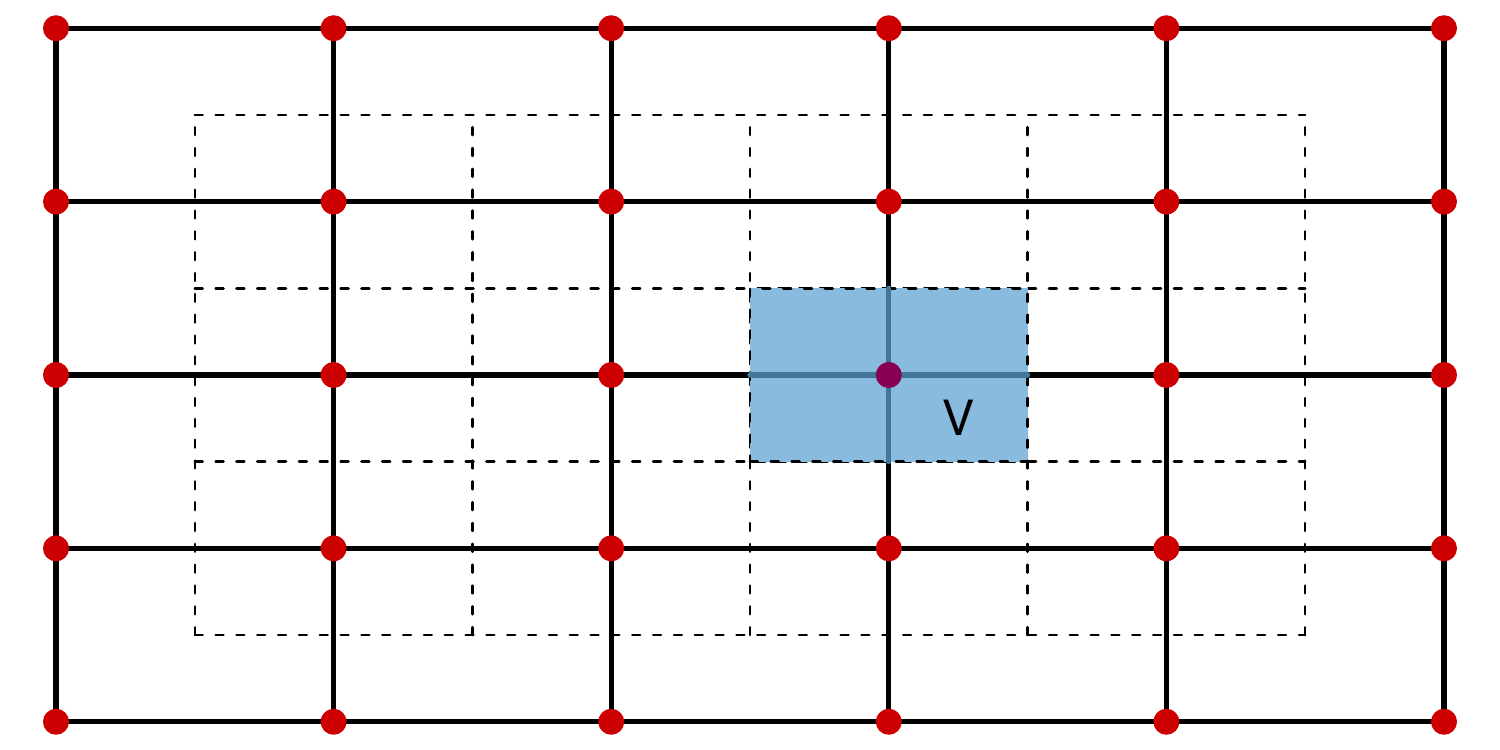}
\caption{Schematic representation of finite volume discretization of a rectangular domain by regular cells
for each inner node.}
\label{fig:finitevols}
\end{figure}

 Define
\[
\mathbf{s}=s\mathbf{I}.
\]
Integrating Biot's equations in a finite volume $V$ (Figure \ref{fig:finitevols}) and considering the terms involving $\sigma$ and $s$, we have

\[
\int_V \nabla\cdot\sigma = \int_{\partial V}\sigma\mathbf{n}
\]
and
\[
\int_V \nabla s = \int_V \nabla\cdot\mathbf{s}=\int_{\partial V}\mathbf{s}\mathbf{n}.
\]

Hence, in the part of the boundary of $V$  contained in the interface, the boundary conditions (\ref{eq:boundary_conds_new}) are straightforward.

\bigskip

\noindent\textbf{Synthetic example. } We consider a water saturated porous medium, also immersed in water. The bone specimen is $4mm$ thick and $10mm$ long. The physical parameters of water are $\rho= 1000 \, Kg/{m^3}$, $K= 2.2\times10^9 \, Pa$.
The source term is as in Nguyen \& Naili \cite{NguyenNaili}, namely

\begin{equation}
F(t)=F_0 e^{-4 \left( f_c t-1 \right)^2 } \sin\left( 2\pi f_c t \right),
\end{equation}
where $f_c=1\, MHz$ and $F_0=1\, m/{s^2}$. The transmiter is  $2\mathrm{mm}$ above the specimen opposite to the receiver. For all experiments a time interval of length  $T = 7\times 10^{-5}s$ is used.

Physical parameters for the porous medium are given in Table 1. 

\begin{table}
\centering
\begin{tabular}{c@{\qquad}c}  \hline
\multicolumn{1}{c}{\bf Parameter} & \multicolumn{1}{c}{\bf Value} \\ \hline
Porosity ($ \phi$)              & $0.5$             \\ 
Tortuosity ($\alpha$)           & $1.4$            \\
Solid bulk modulus $(K_s)$      & $ 20\times 10^9\;\mathrm{Pa}$     \\
Squeletal frame bulk modulus $(K_b)$      & $ 3.3\times 10^9\;\mathrm{Pa}$    \\
Shear modulus $(N)$      & $ 2.6\times 10^9\;\mathrm{Pa}$      \\
Solid density $(\rho_s)$  &  $1960 \; \mathrm{Kg}/\mathrm{m}^3 $
     \\ \hline                         
\end{tabular}
\caption{Physical parameters of porous medium.}\label{tab:parametro}
\end{table}

\begin{figure}
\centering
  \begin{subfigure}
    \centering
    \includegraphics[width=0.45\textwidth]{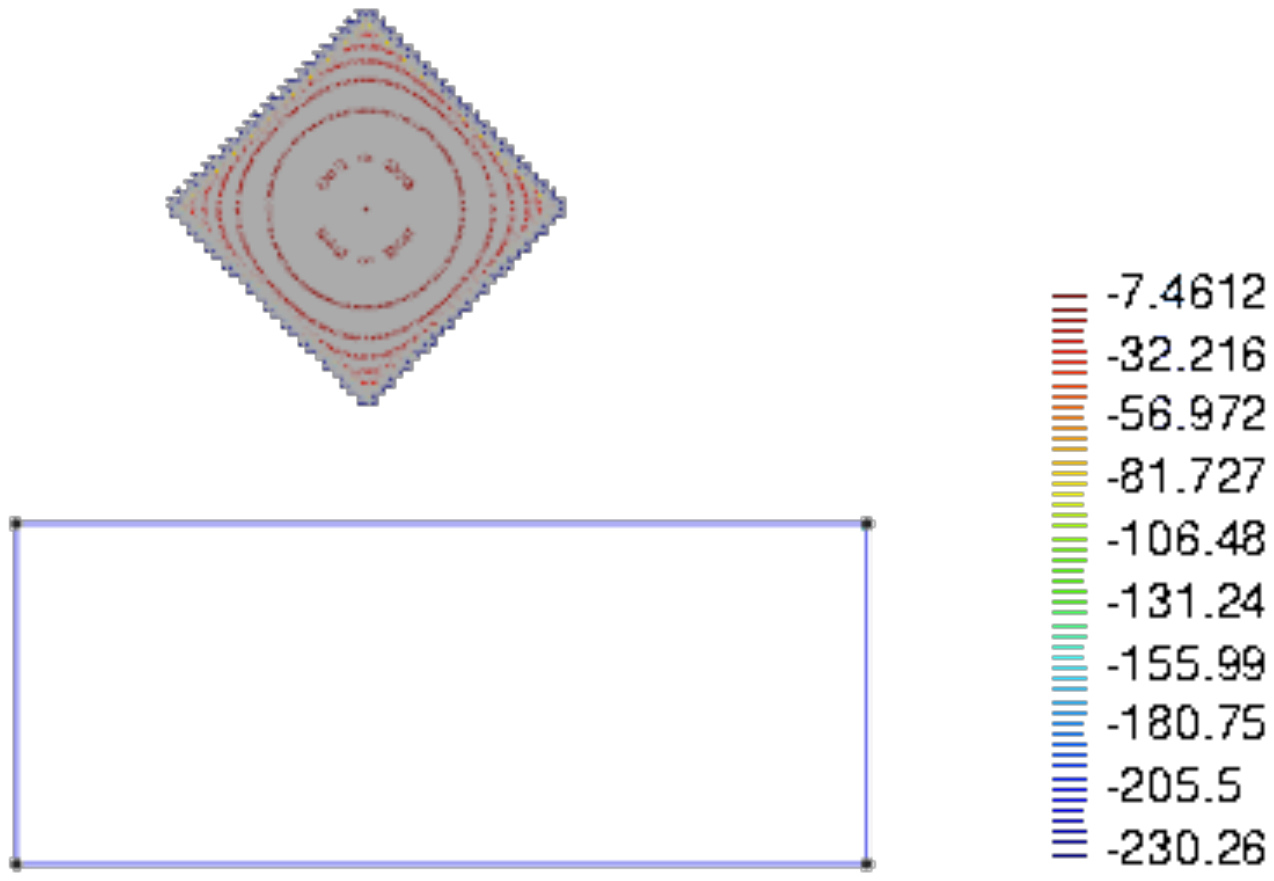}
  \end{subfigure}%
  \begin{subfigure}
    \centering
    \includegraphics[width=0.45\textwidth]{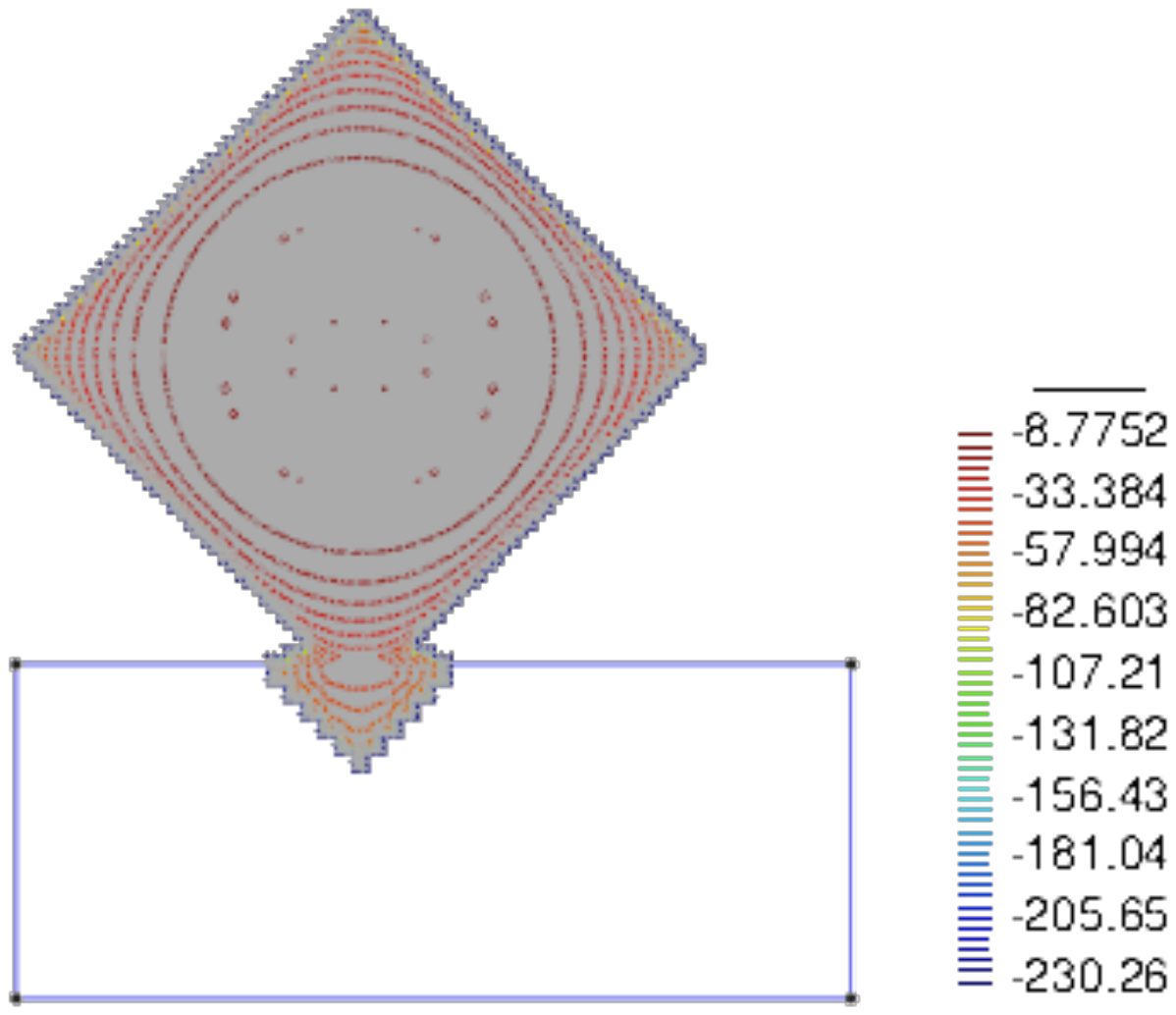}
  \end{subfigure}%
 
  \begin{subfigure}
    \centering
    \includegraphics[width=0.459\textwidth]{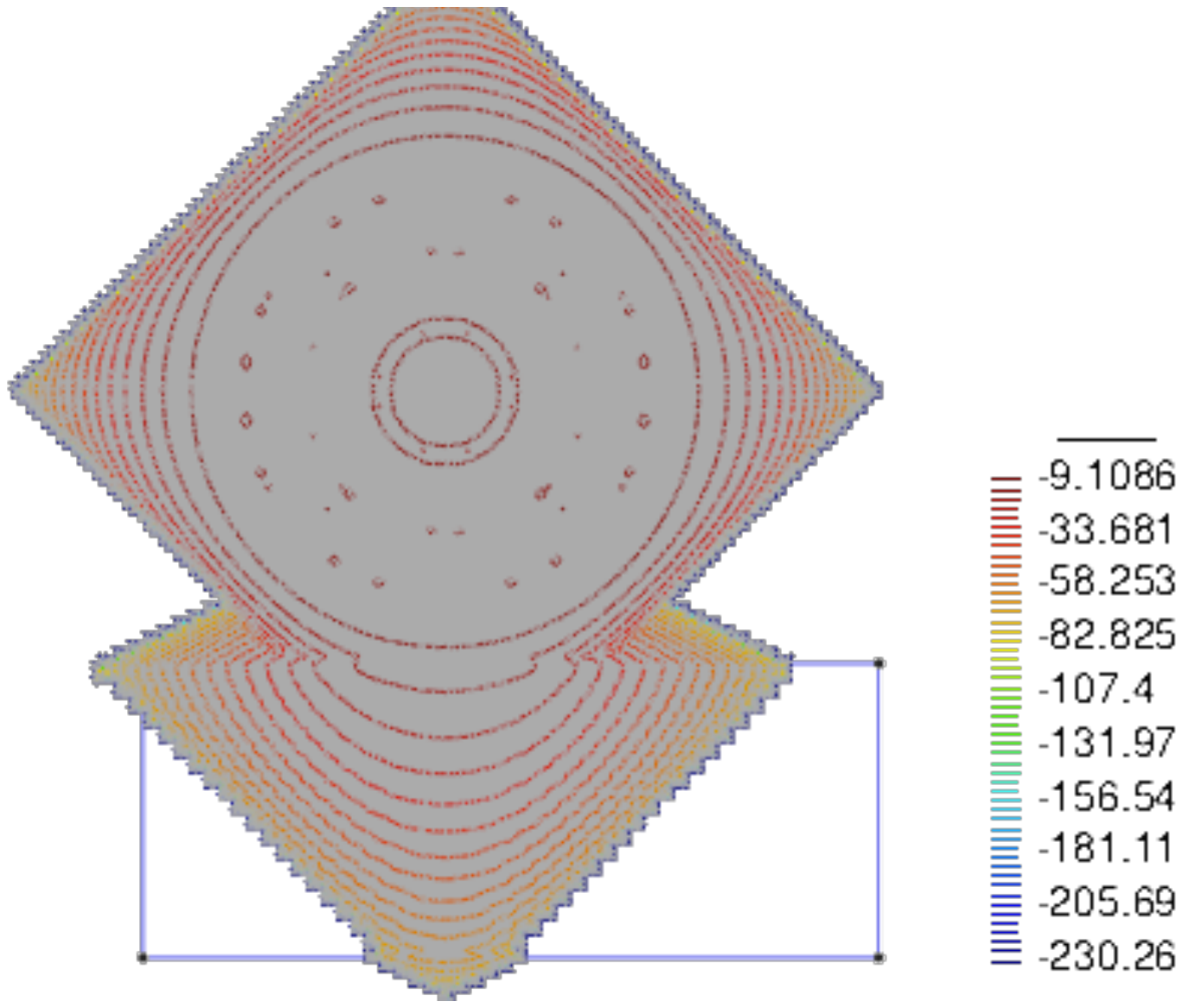}
  \end{subfigure}%
  \begin{subfigure}
    \centering
    \includegraphics[width=0.45\textwidth]{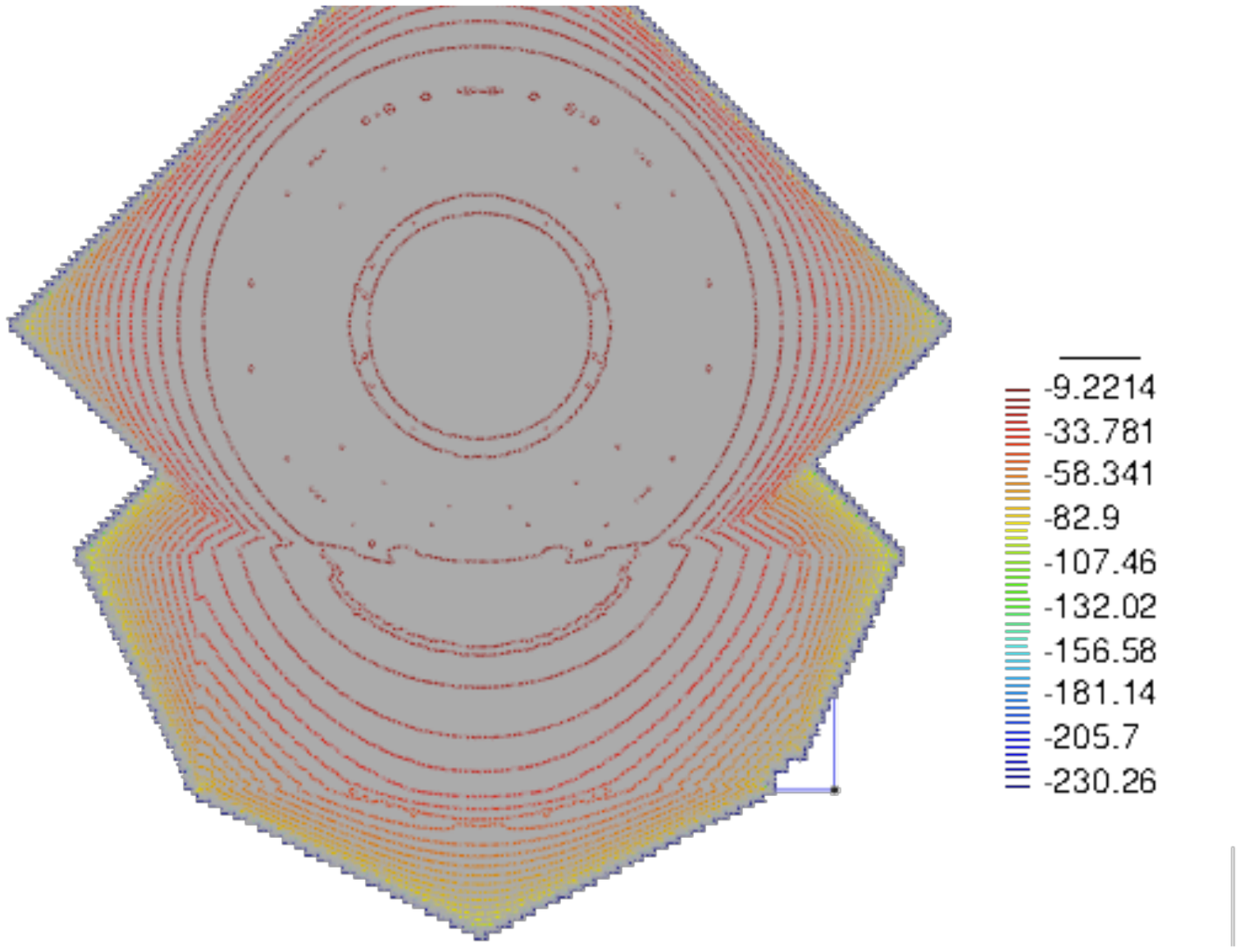}
  \end{subfigure}%
  
  \caption{Evolution of pressure field.}\label{fig:1_uno_1}
\end{figure}

In Figure \ref{fig:1_uno_1} there are some snapshots of the numerical solution obtained by the finite volume method.

\bigskip

\noindent\textbf{Remark. }Numerical modeling of ultrasound propagation trough a heterogeneous, anisotropic, porous material, such as a bone, is a research problem in itself. There is a vast literature on a variety of models and numerical methods. Finite Element in time and frequency domain have been used, as well as Finite Differences. For instance, see Nguyen \& Naili \cite{NguyenNaili} for a hybrid Spectral/FEM approach, Nguyen, Naili \& Sansalone \cite{NguyenNailiSansa} for FEM in time domain, and Chiavassa \& Lombard \cite{ChiavassaLombard} for a Finite Difference implementation. The model above is somewhat simple and does not include all the mechanical complexities. Nevertheless, it is realistic enough to test our methodology of parameter estimation.  It will become apparent that it is straightforward to replace the underlying forward map with more complex models.

\subsection{The parameter estimation problem}

For practical motivations, most studies are conducted on water-saturated specimens rather than medulla, the actual fluid saturating trabecular bone. Consequently, it is customary as a first approximation to consider the fluid saturating the porous medium in Biot's model as known. 

\bigskip

We are led to the inverse problem: Given pressure data
\[
P_i \sim P(x^s,y^s,t_i),\quad i=1,2,\ldots,m,
\]
determine the Biot's parameters
\begin{equation}
\mathbf{u}=\left( \phi\text{, }  \alpha\text{, } K_s\text{, }  K_b\text{, } N\text{, } \rho_s \right).
\end{equation} 

\bigskip

\noindent\textbf{Remark. } It is critical to have realistic reference values in Biot's model, regardless of the chosen methodology for estimation. Obtaining data from lab measurements is complex and costly. For instance, porosity can be obtain from 3D microtomography ($\mu$CT) Wear et al \cite{wearparam}, Pakula et al \cite{pakulaparam}. Tipical values for human trabecular bone vary between $0.55$ and $0.95$ depending on anatomic location and bone situation. Tortuosity values are also scarce. It can be measured using electric spectroscopy, wave reflectometry or estimated from porosity, Hosokawa \& Otani  \cite{hosoparam,hosoparam2}. Reported values are in the interval $[1.01\text{, } \ 1.5]$ Laugier \& Haeiat \cite{bqu}.
Elastic properties of bone tissue are required to estimate macroscopic elastic properties of the saturated skeletal frame. These have been measured using atom strength microscopy, nanoidentation or acoustic microscopy. Then, micro mechanic models can be used to determine volume and shear modulus of the solid.

\section{Bayesian parameter estimation in cancellous bone}

In the sense of Hadamard, a problem is well posed, if existence, uniqueness and continuity with respect to data (stability), can be established. For instance, in differential equations, continuity with respect to initial and/or boundary conditions. A problem is ill-posed if any of the conditions fails. 

Classical well posed problems for differential equations are commonly referred as direct problems, in our case, the Initial-Boundary Value Problem for Biot's model (\ref{eqsBiot_perd_1}), (\ref{eq:scalar_wave}), (\ref{P_Bd}), (\ref{eq:boundary_conds_new}), (\ref{P_Ini}), (\ref{U_Ini}).

In practice, one is interested in a property of a system, to be determined form indirect information. The ill-posed Biot's parameters estimation problem is of this sort, and can be regarded as an inverse problem.

\bigskip

In general, a quantity $\mathbf{y} \in \mathbb{R}^{m}$ is measured to obtain information about another quantity $\mathbf{u} \in \mathbb{R}^{n}$. A model is constructed that relates these quantities. The data $\mathbf{y}$ is usually corrupted by noise. Consequently, the inverse problem can be written as,

\begin{equation} \label{eq39}
 \mathbf{y} = f(\mathbf{u},\mathbf{e})
\end{equation}
where $f:\mathbb{R}^{n}\times\mathbb{R}^{k}\rightarrow\mathbb{R}^{m}$ is a function of the model and  $\mathbf{e} \in \mathbb{R}^{k}$ is the noise vector.

\subsection{Bayesian framework}

Let us  develop the Bayesian methodology for statistical inversion. This paragraph is deliberately terse, for details see Kaipio \& Sommersalo \cite{kaipio} and Stuart \cite{stuart}.

The aim of statistical inversion is to extract information on $\mathbf{u}$, and quantify the uncertainty from the knowledge of $\mathbf{y}$ and the underlying model. It is based on the principles:

\begin{enumerate}
    \item All variables are regarded as \emph{random variables}
 \item Information is on realizations  
 \item This information is coded in probability distributions 
 \item The solution to the inverse problem is the posterior probability distribution
\end{enumerate}

As customary in statistical notation, random variables are capital letters, thus \eqref{eq39}  reads 

\begin{equation} \label{eq40}
 Y = f(U,E).
\end{equation}
In this context, the data $y$ is a realization of $Y$. 

In Bayesian estimation all we know about $U$ it is encompassed in a probability density function, the \emph{prior}, $\pi_0(u)$. The conditional probability density function $\pi(y\vert u)$ is the \emph{likelihood function}, whereas the conditional probability density function $\pi^y(u)\equiv\pi(u\vert y)$ is the \emph{posterior}. All densities are related by the Bayes' formula

 \begin{equation}\label{eq38}
  \pi^y(u)= \frac{ \pi_{pr}(u)\pi(y\vert u)}{\pi(y)}.
 \end{equation}
 
Summarizing, solving a inverse problem in the Bayesian framework, consist on the following: 

\begin{enumerate}
 \item With all available information on $X$, propose a prior $\pi_0(u)$. This is essentially a modeling problem
 
 \item Find the likelihood $\pi(y\vert u)$
 
 \item Develop methods to explore the posterior
\end{enumerate}

For the last step we use a Markov Chain Monte Carlo (MCMC) method. More precisely, we use emcee, an affine invariant  MCMC ensemble sampler. Foreman-Mackey et al \cite{Foremanetal}.

It is well known that this sampling methodology does not depend on the normalizing constant $\pi(y)$ and we write,

\begin{equation*}
  \pi^y(u) \propto  \pi_{pr}(u)\pi(y\vert u).
\end{equation*}

\bigskip

\noindent\textbf{Point estimators}

\bigskip

Given the posterior, a suitable value for the unknown variable $U$ is needed. That is, a point estimator.

A maximizer of the posterior distribution is called a \emph{Maximum A Posteriori estimator}, or MAP estimator. This amounts to the solution of a global optimization problem.

\begin{equation} \label{eq32}
    u_{\,_{MAP}} = \arg \max_{u \in \mathbb{R}^{n}}\pi(ux\vert y),
\end{equation}

Also of interest is the \emph{conditional mean} (MC) estimator, namely

\begin{equation} \label{eq33}
    u_{\,_{MC}} = E(u | y) = \int_{\mathbb{R}^{n}}u \pi(u | y) dx,
\end{equation}

This is an integration problem usually in high dimensions, computationally costly. Classical quadrature rules are prohibited. 

\bigskip

\noindent\textbf{From MAP to Tikhonov}

\bigskip

For Bayesian estimation we consider $\mathbf{u}$ as a random vector distributed as $\pi_0$, a given prior density.
Thus $\vec{y}$ is given by

\[
\vec{y}=\mathcal{G}(\vec{u})+\eta
\]
where $\eta$ is random noise with density $\rho$. Here $\mathcal{G}$ is the observation operator.

From Bayes' formula, the posterior distribution $\pi^y(\vec{u})$ satisfies

\[
    \pi^y(\vec{u}) \propto \rho(\vec{y}-\mathcal{G}(\vec{u}))\pi_0(\vec{u}).
\]

\bigskip

Assuming Gaussian prior $\vec{u}\sim \mathcal{N}(\vec{u}_0,\sigma^2\mathbf{I})$, and Gaussian noise
$\eta\sim \mathcal{N}(\vec{0},\gamma^2\mathbf{I})$ we have

\[
\pi_0(\vec{u})\propto \exp\left( - \frac{1}{\sigma^2}\| \vec{u}-\vec{u}_0\|^2 \right)
\]
and
\[
\rho(\mathbf{u})\propto \exp\left( - \frac{1}{\sigma^2}\| \vec{u}\|^2 \right).
\]
We are led to

\begin{eqnarray*}
    \mathbf{u}_{\,_{MAP}}
    &=& \arg \max\left[\;\exp\left( - \frac{1}{\gamma^2}\| \vec{y}-\mathcal{G}(\vec{u})\|^2 \right)  
\exp\left( - \frac{1}{\sigma^2}\| \mathbf{u}-\vec{u}_0\|^2 \right)\;\right]
\\[5pt]
&=& \arg \min\left[\;\|\vec{y}-\mathcal{G}(\vec{u})\|^2 +
\left(\frac{\gamma}{\sigma}\right)^2\| \mathbf{u}-\mathbf{u}_0\|^2 \;\right].
\end{eqnarray*}

Consequently, $\mathbf{u}_{\,_{MAP}}$ coincides with Tikhonov' solution with regularization parameter $\alpha=\left(\frac{\gamma}{\sigma}\right)^2$.

\subsection{Application to Biot's problem}

\begin{figure}
\centering
\includegraphics[width=160pt]{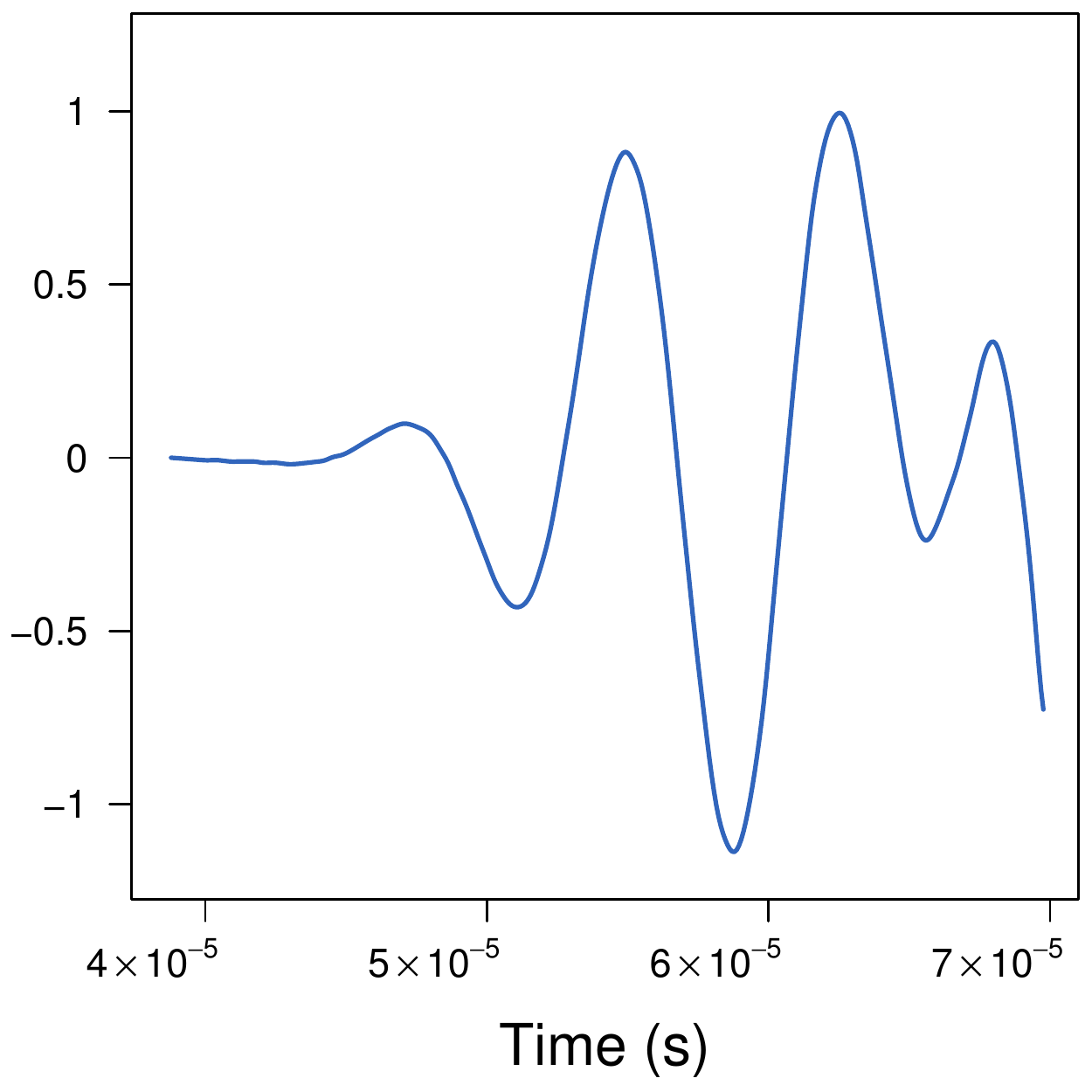} \;
\includegraphics[width=160pt]{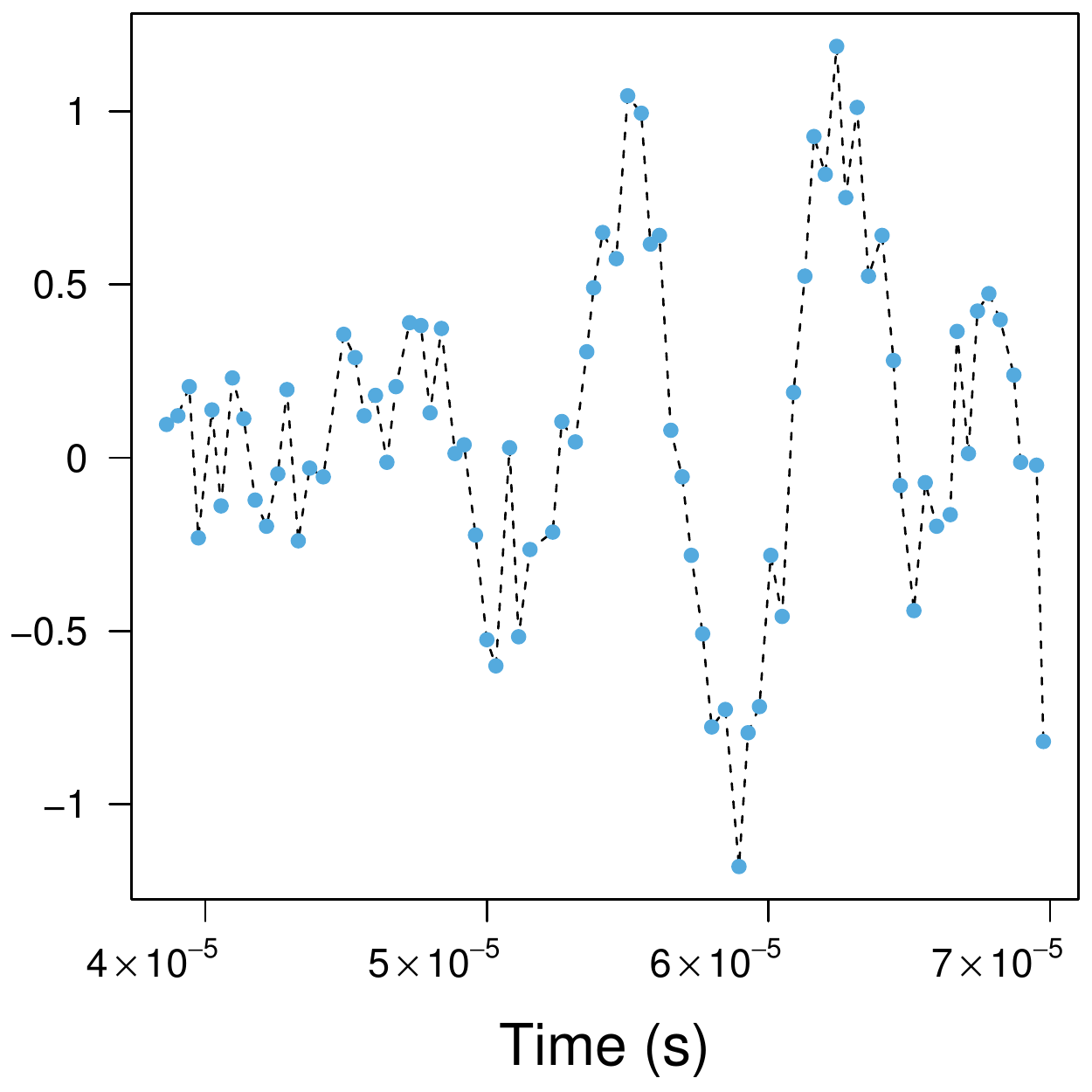}
\caption{Left: Solution of the Biot's model on the position of the receiver.
    Right: Input data are generated by sampling the receiver signal and adding Gaussian noise.}
\label{fig:sintetica_1_noisy}
\end{figure}

We consider the problem described in Section 2.4 and assume  as input data a signal corrupted by Gaussian noise as shown in Figure \ref{fig:sintetica_1_noisy} .

\noindent\textbf{Likelihood function}

\bigskip

In our synthetic example we have considered additive noise, a common assumption. Thus the model becomes,

\begin{equation}\label{eq35}
 Y = f(U) + E.
\end{equation}

Also, it is assumed that the random variables $U$ and $E$ are independent. Consequently, the likelihood function 
is

\begin{equation}\label{eq36}
 \pi(y|u) = \rho(y - f(u)),
\end{equation}
where $\rho$ is a Gaussian probability density function of $E$.

We are led to explore the posterior,

$$\pi(u | y)\propto \pi_{pr}(u)\rho(y - f(u)).$$

\bigskip

\noindent\textbf{Prior densities}

\bigskip

We consider two cases, Gaussian and uniform (uninformative) priors.  In both cases the posterior is sampled and the conditional mean estimator is computed. A comparison is made with the MAP estimator. 

\bigskip

\textbf{Gaussian prior}

\bigskip

In Table \ref{tab:objetivo_gauss_tab} we list the parameters for the a priori Gaussian densities. Notice that the mean of the Gaussian distribution for each parameter is \emph{far away} form the true value. 
 
\begin{table}
\centering
\begin{tabular}{c@{\quad}c@{\quad}c@{\quad}c}
\hline
Property & True value & Mean ($u^0_i$) & Standard deviation ($\gamma_i$) \\ \hline
Porosity ($ \phi$)           & $0.50$             & $0.8$                                           & $ 0.10$  \\ 
Tortuosity ($\alpha$)        & $1.4$        & $ 1.6$                                           & $ 1.5$      \\
Solid bulk modulus $(K_s)$   & $20\times10^9\,$ Pa    & $ 25\times10^9\,$ Pa     & $ 9\times10^9\,$ Pa    \\
Squeletal frame bulk modulus $(K_b)$      & $ 3.3\times10^9\,$ Pa   & $ 3.8\times10^9\,$ Pa     & $ 2.5\times10^9\, $ Pa    \\
Shear modulus $(N)$       & $2.6\times10^9\, $ Pa  & $ 4.5\times10^9\, $ Pa     & $ 5.5\times10^9\, $ Pa    \\
Solid density $(\rho_s)$  & $1960 \, \text{Kg/m}^3 $ &  $1940 \, \text{Kg/m}^3 $  &  $250  \, \text{Kg/m}^3$   \\ \hline                         
\end{tabular}
\caption{Parameters for Gaussian priors.}\label{tab:objetivo_gauss_tab}
\end{table}
  
\begin{figure}
\centering
  \subfigure{
    \centering
   \includegraphics[width=0.48\textwidth]{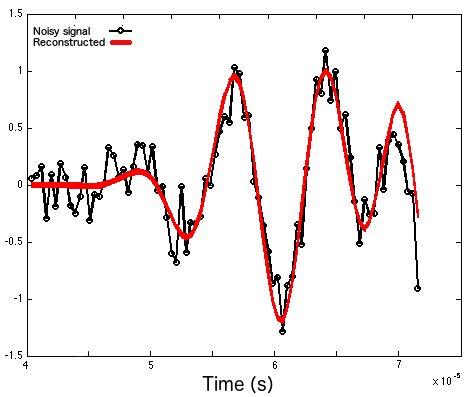}
}%
   \subfigure{
    \centering
    \includegraphics[width=0.48\textwidth]{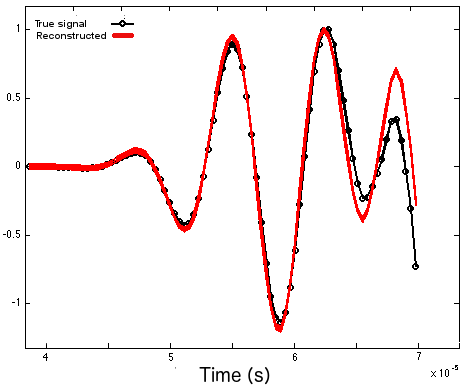}
}%
  \caption{Left: Comparison between the generated signal (solid line) using $u_{\mathrm{CM}}$  and the noisy input (dashed line); Right: Comparsion between the true signal (dashed line) and the recovered signal (solid line).}
  \label{fig:optimized_walks_1_gauss}
\end{figure}

It is remarkable that the conditional mean estimator is capable of recovering the noisy signal. See Figure \ref{fig:optimized_walks_1_gauss}.

\bigskip

\textbf{Uniform prior}

\bigskip

Let us consider uniform priors, namely 
\begin{equation}
\pi(u) \;\propto\; \chi_{[a,b]}(u) = 
\begin{cases} 
1 &\mbox{if } a\leq u\leq b, \\ 
0 & \text{elsewhere}, \end{cases}
\end{equation}
where the parameter of interest is believed to belong to the interval $[a\text{,}b]$. 

The intervals are chosen to be physically meaningful, see Table 3.

\begin{table}
\centering
\label{tabla:propsU}
\begin{tabular}{c@{\qquad}c}
\hline
\multicolumn{1}{c}{\bf Property} & \multicolumn{1}{c}{\bf Interval}  \\ \hline
Porosity ($ \phi$)                                       & $[0.3,\; 0.95]$ \\ 
Tortuosity ($\alpha$)            & $[1,\;\infty)$           \\
Solid bulk modulus $(K_s)$      & $ [1.5 \times 10^{10},\; 3.0 \times 10^{10}] \, $ Pa     \\
Squeletal frame bulk modulus $(K_b)$      & $[2.0 \times 10^{9},\; 4.5\times 10^{9}] \, $ Pa         \\
Shear modulus $(N)$      & $ [2.0\times 10^{9},\; 3\times 10^{9}]\, $  Pa      \\
Solid density $(\rho_s)$  &  $[1000,\; 3000]\; \text{Kg/m}^3 $    \\ \hline                         
\end{tabular}
\caption{Parameters for uniform priors.}\label{tab:objetivo_uniform_tab}
\end{table}

Again, the conditional mean estimator fits satisfactorily even the noisy signal. See Figure 6.

\begin{figure}
\centering
  \subfigure{
    \centering
   \includegraphics[width=0.48\textwidth]{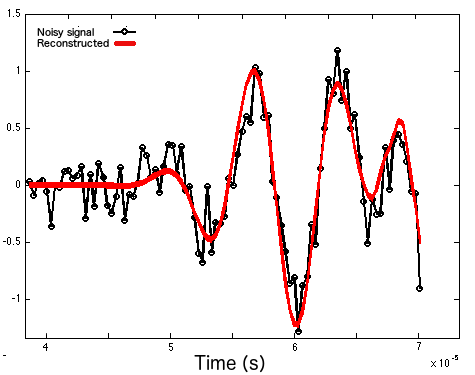}
  }%
   \subfigure{
    \centering
    \includegraphics[width=0.48\textwidth]{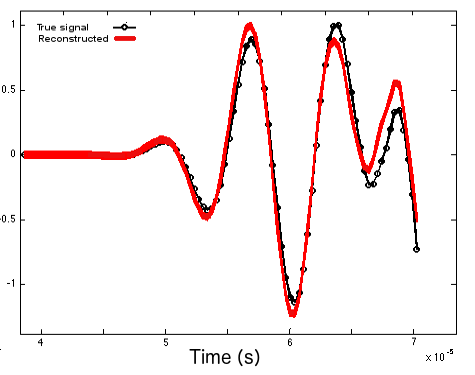}
 }%
 \caption{Left: Comparsion between the reconstructed signal using $u_{\mathrm{CM}}$ (solid line) and the noisy input (dashed line); Right: Comparsion between the reconstructed signal (solid line) and the true signal (dashed line).}
 \label{fig:optimized_walks_uniformative}
\end{figure}

\bigskip

\noindent\textbf{Confidence intervals}

\bigskip 

Having the posterior density, allows to quantify the uncertainty of the estimated parameters. Here we just provide confidence intervals with a 0.9 probability. Results are shown in Tables 4 and 5 for Gaussian and uniform priors respectively.

\begin{table}
\centering
\begin{tabular}{@{\quad}l@{\quad}llll@{\quad}l@{\quad}} \hline
   &   \multicolumn{3}{c}{Gaussian prior density }    &  & \\
Parameter   & $u_{\mathrm{TRUE}}$ &   & $u_{\mathrm{CM}}$ & Interval  \\ \hline
$\phi$    & $0.5$ &  &     $0.536 $   &   $[\  0.478 \text{,}\  0.614 \  ]$        \\
$\alpha$ & $1.4$ &       &   $1.421$   &  $[\  1.350  \text{,}\  1.505 \  ]$        \\
$K_s$ ($\times10^ {10}$Pa)        & $2.0$  &      &   $2.000 $   &    $[\ 1.9999922  \text{,}\  2.0000079\  ] $     \\
$K_b$ ($\times10^ {9}$Pa)         & $3.3$  &     &  $3.299$   &    $[\  3.299  \text{,}\  3.300  \  ]$        \\
$N$  ($\times10^ {9}$Pa)          & $2.6$  &     & $2.544$  &      $[\  2.226  \text{,}\  2.862  \  ]$      \\
$\rho_s$  ($\times10^ {3}$Kg/$m ^3$)   & $1.96$  &      &  $1.955$    &    $[\ 1.943\text{,}\  1.969 \  ]$        \\ \hline
\end{tabular}
\caption{Comparison of estimated parameters.}
\label{tab:CM_G}
\end{table}

\begin{table}
\centering
\begin{tabular}{@{\quad}l@{\quad}llll@{\quad}l@{\quad}}
\hline
   &   \multicolumn{3}{c}{Uniform prior density }  &  &  \\
Parameter   & $u_{\mathrm{TRUE}}$ &  &  $u_{\mathrm{CM}}$  & Interval \\ \hline
$\phi$    & $0.5$ &     &    $0.549 $  &  $[\  0.505   \text{,}\  0.642   \  ]$      \\
$\alpha$ & $1.4$ &  &       $1.432$  &  $[\  1.321  \text{,}\  1.540   \  ]$          \\
$K_s$ ($\times10^ {10}$Pa)        & $2.0$  &  &     $2.085$  &    $[\  1.477  \text{,}\  2.451 \  ]$   \\
$K_b$ ($\times10^ {9}$Pa)         & $3.3$  &  &         3.270  &  $[\  2.847    \text{,}\  3.808   \  ]$    \\
$N$  ($\times10^ {9}$Pa)          & $2.6$  &  &       $2.682$  &  $[\  1.645  \text{,}\  3.482   \  ]$    \\
$\rho_s$  ($\times10^ {3}$Kg/$m ^3$)   & $1.96$  &   &    $1.949$ & $[\  1.815   \text{,}\  2.077  \  ]$  \\ 
\hline
\end{tabular}
\caption{Comparison of estimated parameters.}
\label{tab:CM_U}
\end{table}

\bigskip

\noindent\textbf{Remark. }A drawback of Bayesian estimation is its computational cost. En each step of the random walk of the MCMC method, the forward map involves the solution of the Biot's model. Nevertheless, the results above show that the conditional mean is a reliable point estimator. We shall see below that in this case other approaches of estimation may not suffice.

\section{A PDE-Constrained optimization approach}

Let $\pi_0(\vec{u})$ denote the prior density for the Biot's parameters. As pointed out above, a classical approach is to consider the problem of estimation as one of optimization by means of the MAP estimator. Namely, 

\begin{equation}\label{eq:opti_pi}
\vec{u}_{\,_{MAP}} = \arg \max_{\vec{u}} \ \exp \left\lbrace {-\frac{1}{\sigma  }|| y - \mathcal{G}(\vec{u}) ||^2}   \right\rbrace \pi_0(\vec{u}).
\end{equation}

Applying natural logarithm 
\begin{equation}\label{eq:opti_5}
\vec{u}_{\,_\mathrm{MAP}} =  \arg \min_{\vec{u}} \  || y - \mathcal{G}(\vec{u}) ||^2 - \sigma \log\left( \pi_0(\vec{u})\right).
\end{equation}

This minimization problem is constrained by the IBVP for Biot's model (\ref{eqsBiot_perd_1}), (\ref{eq:scalar_wave}), (\ref{P_Bd}), (\ref{eq:boundary_conds_new}), (\ref{P_Ini}), (\ref{U_Ini}).  Consequently any evaluation of the objective function requires the solution of this IBVP.
Derivative based algorithms are computationally expensive. A reasonable alternative is the Nelder-Mead method, a description below.

\subsection{The Nelder--Mead method}

The optimization problem to calculate the MAP estimator is solved by using the derivative free Nelder--Mead method
\cite{nelderMead}.

The method only requires evaluations of the objective function in \eqref{eq:opti_5}.
It is based on the iterative update of a \emph{simplex}, which in this case
is a set of $n+1$ points in  $n$-dimensional space and they do not lie in a space of lower
dimension. Each point of the simplex is called a \emph{vertex}. 
The shape and size of the simplex is modified according to the values of the objective function
$f:\mathbb{R}^n \rightarrow \mathbb{R}$ at each vertex. 

\begin{figure}
\centering
    \includegraphics[width=\textwidth]{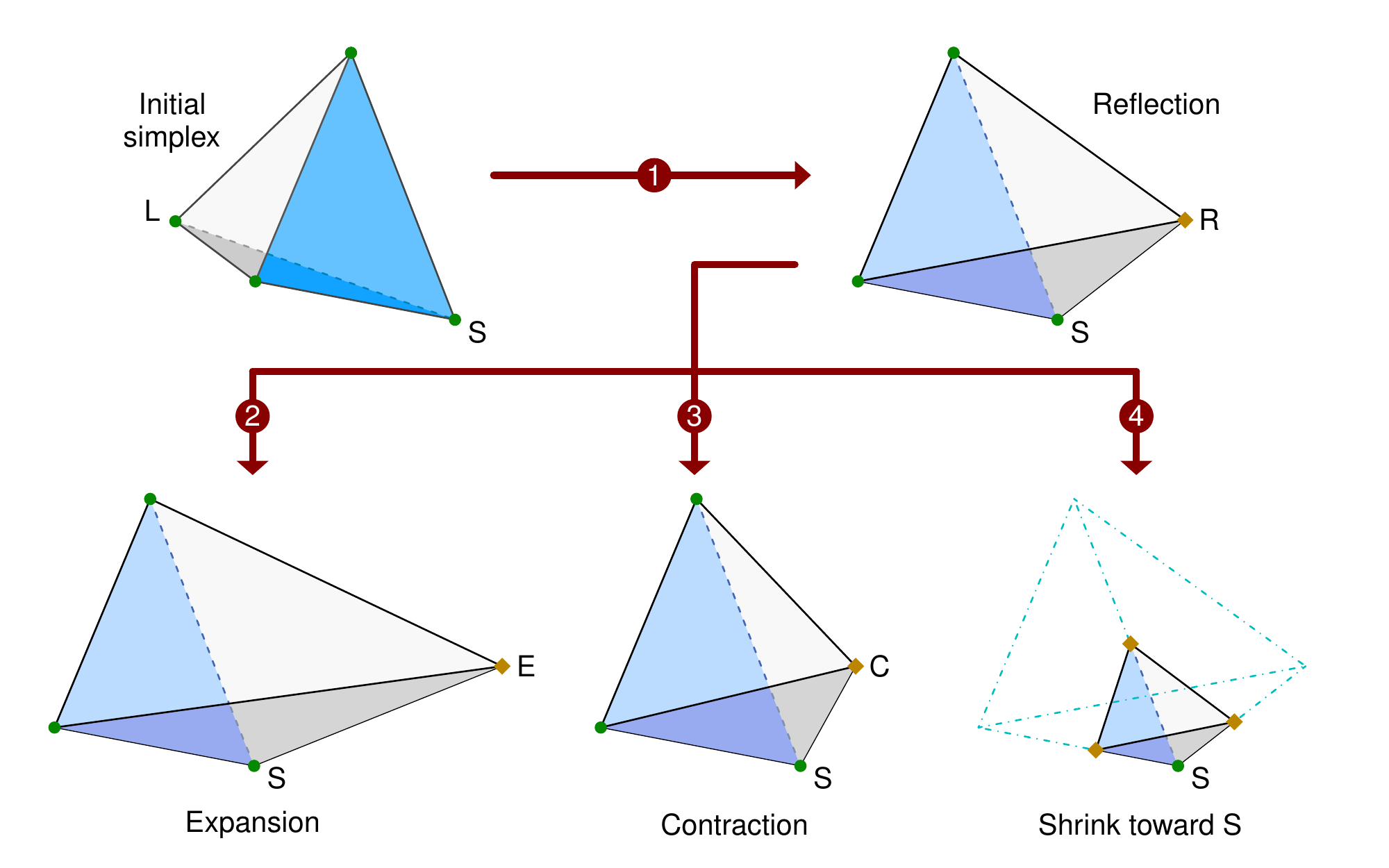}
\caption{Steps of the  Nelder--Mead method.}
\label{fig:nelderMead}
\end{figure}

The algorithm starts with an initial guess of the vertices.In reference to Fig. \ref{fig:nelderMead},
let $L$ and $S$  be the vertices where  the objective function has 
its largest and smallest values, respectively. 
The algorithm tries to modify the vertex $L$ to find 
a new point, such that the value of the objective function at this point will be 
smallest than $f(S)$. To illustrate the complexity of the method we delve a little further.

The Nelder-Mead algorithm has 
four parameters:
\begin{itemize}
\item the reflection coefficient $\tau_r>0$ (usually $\tau_r$ is set to $1$),
\item the expansion factor $\tau_e>\max\{1, \tau_r\}$,
\item the contraction parameter  $\tau_c \in (0, 1)$, and
\item the shrinkage factor $\tau_s\in (0, 1)$.
\end{itemize}

In each iteration all the vertices are indexed according to the values of the function $f(x)$.
Thus the simplex is composed by the vertices 
\begin{equation}
    S=v_1, v_2, ..., v_{n+1}=L \quad \text{if} \quad 
    f_i=f(v_i) \leq f_{i+1}=f(v_{i+1}), \quad  i=1,2, ..., n.
    \label{eqn:indexing}
\end{equation}

To define the transformations of the simplex, we need to calculate
the centroid $\bar{v}$ of the $n$ first vertices,
      \begin{equation}
        \bar{v} = \frac{1}{n} \sum_{i=1}^n v_i,
        \label{eqn:centroid}
      \end{equation}
and the point 
\[ v(\tau) = \bar{v} + \tau\,(\bar{v} - S), \]
for  $\tau \in \{\tau_r, \tau_e, \tau_c, \tau_s \}$. 
Each parameter is associated to an operation (see Fig. \ref{fig:nelderMead}):
\begin{itemize}
\item \emph{Reflection} produces a movement of the simplex towards regions where $f$ is getting smaller values.
\item \emph{Expansion} increases the size of the simplex to advance more quickly in search of the local minimum.
\item \emph{Contraction} is applied when reflection and expansion fails, and it allows to get an inner point of 
the simplex in which $f$ takes a value lower than $f(L)$ at least. 
\item \emph{Shrink toward $S$} moves all the vertices in the direction of the current best point to reduce the size of the simplex when it is in a valley of the objective function.  This allows previous operations can continue to be applied in the following iterations.
\end{itemize}

The algorithm applies the following steps in each iteration to find a local minimum
of the function $f(x)$: 

\begin{enumerate}
    \item Indexing the vertices of the simplex  according to the objective function values \eqref{eqn:indexing}.

\item Calculate the centroid $\bar{v}$ of the first $n$ vertices 
      \eqref{eqn:centroid}.
     
\item Transform the simplex by the following operations:

\begin{enumerate}
\item (Reflection) Calculate  $R = v(\tau_r)$. If $f(S) < f(R) < f(L)$, $L$ is replaced by $R$
      and we move to the step 1 to start the next iteration.

\item (Expansion) If $f(R) < f(S)$,  we calculate  $E = v(\tau_e)$.
      If $f(E)<f(S)$, $L$ is replaced by $E$.  
      Otherwise  $L$ is replaced with $R$. The process is restarted. 
      This 

\item (Contraction) If $f(v_i) < f(R)$ for $i=1,...,n$, the point 
      $C = v(\tau_c)$
      is calculated. \\
      If $f(C) < \min\{f(L), f(R) \}$, $L$ is replaced with $C$ and the process is restarted.
      Otherwise a shrink toward $S$ is applied.

\item (Shrink toward $S$) For $i=2,3,...,n+1$, the vertices $v_i$ are modified by
      \[ v_i = S + \tau_s \,(v_i - S).   \]

\end{enumerate}

\end{enumerate}

The process continues until a maximum number of iterations is reached, or
when the simplex reaches some minimum size, or the best vertex $S$ becoming less
than a given value or failing to change its position between successive iterations.

\bigskip

\begin{figure}
\centering
    \includegraphics[width=0.6\textwidth]{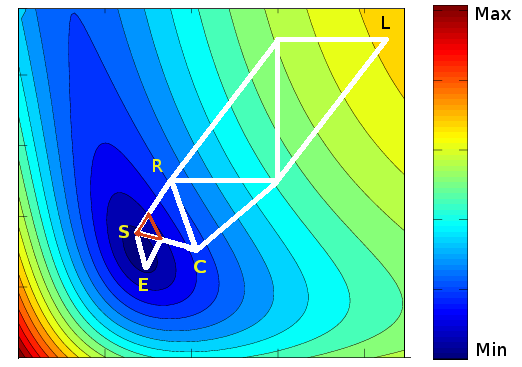}
\caption{Illustration of the evolution of a simplex to find the minimum value of a function
    $f:\mathbb{R}^2 \rightarrow \mathbb{R}$.}
\label{fig:nelderMeadSteps}
\end{figure}

Fig. \ref{fig:nelderMeadSteps} shows
a 2D example of the transformations applied to an initial simplex. 
First  a reflection of the $L$ vertex is applied to reach the point $R$. Next 
a contraction replaces the vertex $R$  by the  point $C$. Then a contraction followed by a shrink contraction 
reduces the size of simplex to allow it reaches a valley floor and finally an expansion is applied.

\noindent\textbf{Remark. }(i) We contend that the Nelder-Mead method is appropriate for the problem at hand. 
It does not require the objective function to be smooth, hence it does not require computation of derivatives. \newline
(ii) On the down side, it is well known that its performance decreases significantly in problems with more
than $10$ variables Han \& Neumann \cite{hanNeumann}. Also,  in the case of problems with 
few variables it may fail to converge to a critical point of the objective function Mckinnon \cite{mckinnon}.

\subsection{Nelder--Mead solutions to Biot's problem}

Chronologically, we posed the estimation problem as one of PDE-Constrained optimization. In the Biot's problem  the number of parameters to be calculated using the MAP estimator \eqref{eq:opti_5} is at most six. Thus Nelder-Mead is appropriate.  The prior information of the variables was used to build the initial simplex. 

First assuming a gaussian prior, or equivalently a regularized least square problem, the MAP estimators in Table 6 are obtained. The true and recovered signal are shown in Figure 9.  Starting with estimating the full set of six parameters, it was observed that the method is unable to recover even the noiseless signal. The problem was simplified one parameter at time. For instance, Figure 9(c) show the estimated signal assuming $\rho_s$ known, and so on. 

\begin{table}
\centering
\begin{tabular}{@{\quad}l@{\quad}llll@{\quad}l@{\quad}} \hline
Parameter   & $u_{\mathrm{TRUE}}$ &  & $u_{\mathrm{MAP}}$   \\ \hline
$\phi$    & $0.5$ &  &    $0.541$         \\
$\alpha$ & $1.4$ &  &  $4.171$             \\
$K_s$ ($\times10^ {10}$Pa)        & $2.0$  &  &  $2.369$        \\
$K_b$ ($\times10^ {9}$Pa)         & $3.3$  &  &   $3.241$          \\
$N$  ($\times10^ {9}$Pa)          & $2.6$  &  &    $9.651$         \\
$\rho_s$  ($\times10^ {3}$Kg/$m ^3$)   & $1.96$  &   &   $2.149$          \\ \hline
\end{tabular}
\caption{Comparison of estimated parameters for a Gaussian prior.}
\label{tab:MAP_MC_G}
\end{table}

\begin{figure}
\centering
  \subfigure{
    \centering
    \includegraphics[width=0.48\textwidth]{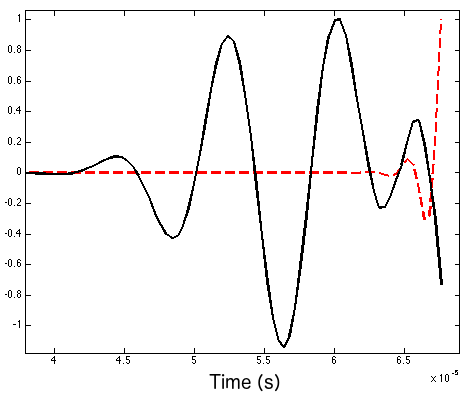}
    \label{fig:ex3-a}
  }%
   \subfigure{
    \centering
    \includegraphics[width=0.48\textwidth]{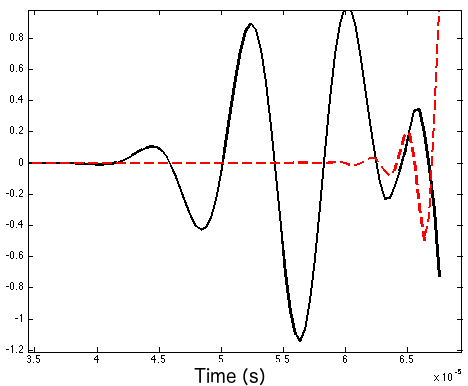}
    \label{fig:ex3-b}
}%
\\
  \subfigure{
    \centering
    \includegraphics[width=0.48\textwidth]{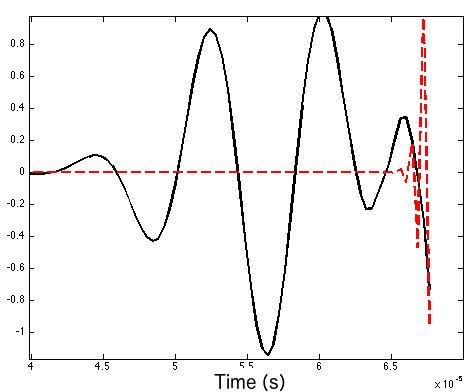}
    \label{fig:ex3-c}
}%
   \subfigure{
    \centering
    \includegraphics[width=0.48\textwidth]{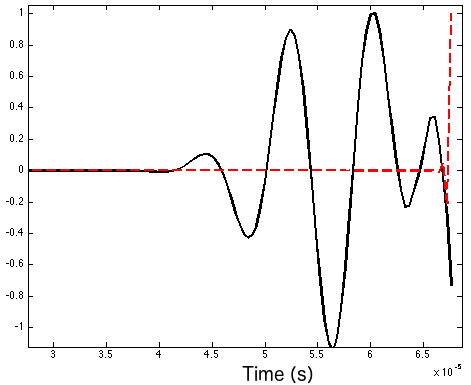}
    \label{fig:ex3-d}
}%
\caption{Comparison between four reconstructed signals using $u_{\,_\mathrm{MAP}}$ with Gaussian priors  (dashed line) and the true signal (solid line) for the estimated parameters: (a) $u=(\phi,\ \alpha,\ K_s)$, (b) $u=(\phi,\ \alpha,\ K_s,\ K_b)$, (c) $u=(\phi,\ \alpha,\ K_s,\ K_b,\ N)$, (d) $u=(\phi,\ \alpha,\ K_s,\ K_b,\ N, \rho_s)$.}
\label{fig:ex3}
\end{figure}

Next uniform priors are considered and the same experiment is carried out. As shown in Table 7 and Figure 10, no improvement is attained. We remark that other optimization methods also fail. 

\begin{table}
\centering
\begin{tabular}{@{\quad}l@{\quad}llll@{\quad}l@{\quad}}
\hline
Parameter   & $u_{\mathrm{TRUE}}$ &  & $u_{\mathrm{MAP}}$  \\ \hline
$\phi$    & $0.5$ & &   $0.613$           \\
$\alpha$ & $1.4$ &  &   $1.350$            \\
$K_s$ ($\times10^ {10}$Pa)        & $2.0$  &  &    $1.382$     \\
$K_b$ ($\times10^ {9}$Pa)         & $3.3$  &  &     $2.301$      \\
$N$  ($\times10^ {9}$Pa)          & $2.6$  &  &     $4.212$        \\
$\rho_s$  ($\times10^ {3}$Kg/$m ^3$)   & $1.96$  &   &  $2.351$      \\ 
\hline
\end{tabular}
\caption{Comparison of estimated parameters for a uniform prior.}
\label{tab:MAP_MC_U}
\end{table}

\begin{figure}
\centering
  \subfigure{
    \centering
    \includegraphics[width=0.48\textwidth]{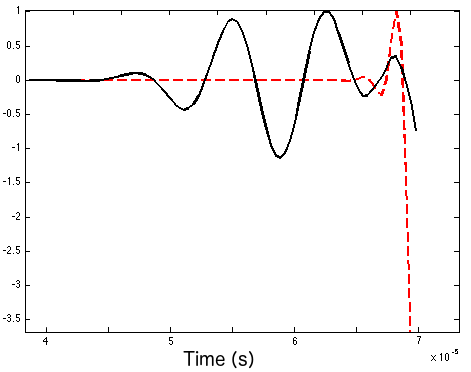}
 }%
   \subfigure{
    \centering
    \includegraphics[width=0.48\textwidth]{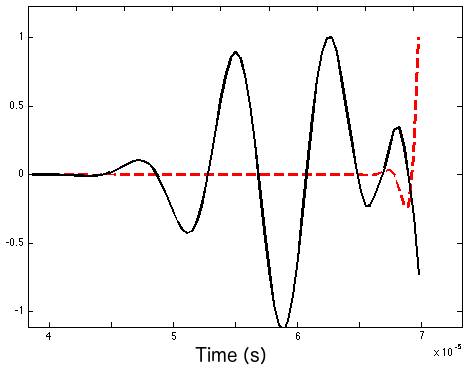}
}%
  
  \subfigure{
    \centering
    \includegraphics[width=0.48\textwidth]{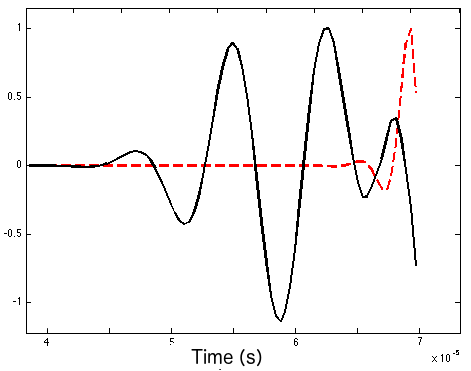}
}%
   \subfigure{
    \centering
    \includegraphics[width=0.48\textwidth]{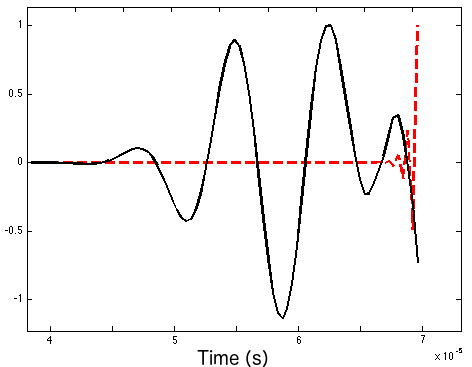}
}%
  
\caption{Comparison between four reconstructed signals using $u_{\mathrm{MAP}}$ with uniform priors  (dashed line) and the true signal (solid line) for the estimated parameters: (a) $u=(\phi,\ \alpha,\ K_s)$, (b) $u=(\phi,\ \alpha,\ K_s,\ K_b)$, (c) $u=(\phi,\ \alpha,\ K_s,\ K_b,\ N)$, (d) $u=(\phi,\ \alpha,\ K_s,\ K_b,\ N, \rho_s)$.}
\label{fig:map_vs_uniform}
\end{figure}

\section{Conclusions}

We have introduced the problem of parameter estimation for a Biot's medium modeling a cancellous bone. The problem has been posed for both a minimization problem and a posterior density estimation in the Bayesian framework. The MAP estimator is shown as the solution of the minimization problem. We carried out extensive experiments with a variety of methods, classical descent methods as well as derivative free methods. All led to the same conclusion, the MAP estimator in not capable of recovering the given signal. We show results only for the derivative free method Nelder Mead. In contrast, the conditional mean recovers the noisy signal satisfactorily. Consequently, although the computation of the conditional mean is costly because of the sampling of the posterior, its use is advisable for diagnostics of the bone properties. 

One may argue that the conditional mean is better suited to represent an intrinsically heterogeneous porous medium, a query worth of an in depth study.  Also of interest, is to consider the initial geophysical phenomena where the Biot's model applies.

\bigskip

\center{Acknowledgements} 

M. A. Moreles would like to acknowledge to ECOS-NORD project number \newline 000000000263116/M15M01 for financial support during this research.

\end{document}